\documentclass[12pt]{article}

\def\be{\begin{equation}}
\def\ee{\end{equation}}
\def\ov{\overline}

\newtheorem{theorem}{Theorem}
\newtheorem{lemma}[theorem]{Lemma}

\newtheorem{definition}[theorem]{Definition}
\newtheorem{proposition}[theorem]{Proposition}

\begin{document}

\title{Ultrametric random field}

\author{A.Yu.Khrennikov\footnote{International Center for Mathematical
Modelling in Physics and Cognitive Sciences, University of
V\"axj\"o, S-35195, Sweden, e--mail:
Andrei.Khrennikov@msi.vxu.se}, S.V.Kozyrev\footnote{Steklov
Mathematical Institute, Moscow, Russia, e--mail:
kozyrev@mi.ras.ru}}

\maketitle

\begin{abstract}
Gaussian random field on general ultrametric space is introduced
as a solution of pseudodifferential stochastic equation.
Covariation of the introduced random field is computed with the
help of wavelet analysis on ultrametric spaces.

Notion of ultrametric Markovianity, which describes independence
of contributions to random field from different ultrametric balls
is introduced. We show that the random field under investigation
satisfies this property.
\end{abstract}

\bigskip

\centerline{Keywords: ultrametric analysis, random fields,
pseudodifferential operators}

\bigskip

\centerline{AMS Subject Classification: 60G15, 60G60, 60H15,
60H40, 60J99}

\section{Introduction}

In the present paper we apply analysis of ultrametric
pseudodifferential operators and ultrametric wavelets to
investigation of Gaussian random fields, which are solutions of
stochastic pseudodifferential equations on general ultrametric
spaces (actually, these spaces obey some properties which
guarantee, for instance, local compactness of these spaces). This
construction is a far generalization of the $p$--adic Brownian
motion of Bikulov and Volovich \cite{BV}, which is a solution of
$p$--adic stochastic pseudodifferential equation
$$
D^{\alpha}\psi(x)=\phi(x)
$$
Here $D^{\alpha}$ is the Vladimirov operator of $p$--adic
fractional derivation \cite{VVZ} and $\phi$ is the white noise on
the field of $p$--adic numbers $Q_p$.

The main result of the present paper is theorem \ref{existence}
(see Section 3), which schematically can be formulated as follows.

\bigskip

\noindent{\bf Theorem}\qquad {\sl Let $X$ be an ultrametric space
(satisfying some natural properties), $\nu$ a measure on $X$, for
which measures of all balls are positive, $T$ be positive
pseudodifferential operator on this space of the form
$$
Tf(x)=\int T{({\rm sup}(x,y))}(f(x)-f(y))d\nu(y)
$$
and let some conditions of convergence be satisfied.

Let $\phi(x)$ be white noise on $X$ (i.e. Gaussian mean zero real
valued $\delta$--correlated generalized random field).

Then the solution of the stochastic pseudodifferential equation
$$
T\psi(x)=\phi(x)
$$
is a Gaussian mean zero random field with covariation
$$
\langle \overline{\psi(x)}\psi(x')\rangle=-\lambda^{-2}_{{\rm
sup}(x,x')}\nu^{-1}({{\rm sup}(x,x')})+\sum_{I> {\rm
sup}(x,x')}\lambda^{-2}_{I}\left(\nu^{-1}({I-1,{\rm
sup}(x,x')})-\nu^{-1}(I)\right)
$$

The corresponding Gaussian measure is $\sigma$--additive on the
$\sigma$--algebra of weak Borel subsets of the space of
distributions $D'(X)$. }

\bigskip

We also introduce the definition of ultrametric Markovianity for
random fields of ultrametric argument, and show that the described
by the theorem above random field is ultrametric Markovian.

In the theorem above ultrametric space $X$ should satisfy some
natural properties (which imply, for instance, local compactness).
The construction of solutions of stochastic pseudodifferential
equations on ultrametric spaces given in this article can be
generalized on arbitrary ultrametric spaces. In the general case
the space $D(X)$ of test functions is not nuclear and one could
not apply the theorem of Minlos--Sazonov. However, even in this
case it is possible to show $\sigma$--additivity of the
corresponding Gaussian measure on $D'(X)$ by using results of the
paper by O.G.Smolyanov and S.V.Fomin \cite{SF}.

The present paper continues development in the field of
ultrametric analysis and $p$--adic mathematical physics, see
\cite{VVZ} for other results in this direction. We use analysis of
general ultrametric pseudodifferential operators and wavelets,
developed in \cite{ACHA}--\cite{nextIzv}.

For exposition of theory of $p$--adic pseudodifferential operators
(PDO) see \cite{VVZ}. It was found \cite{VVZ} that $p$--adic PDO
(which can be diagonalized by the Fourier transform) have also
bases of locally constant eigenvectors with compact support (i.e.
vectors of these bases are functions in the space of test
functions $D(Q_p)$). In \cite{wavelets} an example of such a basis
was found, which is gives the construction of basis of wavelets on
$Q_p$. In \cite{Benedetto}, \cite{Benedetto1} generalizations of
construction of $p$--adic wavelets \cite{wavelets} were built for
locally compact abelian groups. In \cite{nhoper}--\cite{nhoper2}
$p$--adic wavelet bases were used for diagonalization of $p$--adic
PDO, which can not be diagonalized by the Fourier transform. In
\cite{ACHA}--\cite{nextIzv} theory of wavelets and
pseudodifferential operators on general ultrametric spaces (which
do not possess any group structure) was constructed. Introduced
there ultrametric wavelets are eigenvectors of ultrametric PDO.
Therefore, on general ultrametric space there exists an analysis
of PDO and of wavelets, and wavelet transform play the role of the
Fourier transform, used in the analysis on abelian groups.

For other results in ultrametric an $p$--adic mathematical physics
see \cite{Vstring}--\cite{AndrAlb}. For instance, in cite
\cite{Vstring}, \cite{Freund} $p$--adic strings were investigated.
In \cite{ABK}, \cite{PaSu} $p$--adic analysis was applied to
investigation of spin glasses in the replica approach.

The exposition of the present paper is as follows.

In Section 2 we put the exposition of results in ultrametric PDO,
wavelets, and distribution theory.

In Section 3, using the analysis of ultrametric wavelets, we
construct Gaussian random field, which is the solution of
pseudodifferential stochastic equation on general ultrametric
space. We introduce the definition of ultrametric Markovianity and
show that the obtained random field is ultrametric Markovian.

In Section 4 (the Appendix) we discuss the Minlos--Sazonov theorem
about $\sigma$--additivity of measures on topological vector
spaces.

\section{Ultrametric analysis}

In this Section we put the results on ultrametric analysis, which
mainly may be found in \cite{ACHA}, \cite{Izv}, \cite{nextIzv}.

We introduce ultrametric wavelet analysis, analysis of ultrametric
pseudodifferential operators (PDO), and ultrametric distribution
theory.

Rich theory of ultrametric wavelets and pseudodifferential
operators, and especially the fact that ultrametric PDO can be
diagonalized by ultrametric wavelet transform shows that analysis
of ultrametric wavelets plays the role of the Fourier analysis for
general ultrametric spaces. Since the existence of the Fourier
transform is related to structure of abelian group, this shows
that that the ultrametricity property actually is the analogue of
the group property.

\subsection{Ultrametric PDO}

\begin{definition}{\sl
An ultrametric space is a metric space with the ultrametric $|xy|$
(where $|xy|$ is called the distance between $x$ and $y$), i.e.
the function of two variables, satisfying the properties of
positivity and non degeneracy
$$
|xy|\ge 0,\qquad |xy|=0\quad \Longrightarrow\quad x=y;
$$
symmetricity
$$
|xy|=|yx|;
$$
and the strong triangle inequality
$$
|xy|\le{\rm max }(|xz|,|yz|),\qquad \forall z.
$$
}
\end{definition}

We say that ultrametric space $X$ is regular, if this space
satisfies the following properties:

\medskip

1) The set of all the balls of nonzero diameter in $X$ is no more
than countable;

\medskip

2) For any decreasing sequence of balls $\{D^{(k)}\}$,
$D^{(k)}\supset D^{(k+1)}$, diameters of the balls tend to zero;

\medskip

3) Any ball is a finite union of maximal subballs.

\bigskip

For ultrametric space $X$ consider the set ${\cal T}(X)$, which
contains all the balls in $X$ of nonzero diameter, and the balls
of zero diameter which are maximal subbals in balls of nonzero
diameter. Remind that non--maximal vertex $I$ of a directed tree
has the branching index $p_I$, if this vertex is connected to
$p_I+1$ other vertices by links (and to $p_I$ other vertices in
the case of the maximal vertex). The following theorem can be
found in \cite{nextIzv} (the analogous result was obtained in
\cite{Lemin}, see also \cite{Freud} where it was presented).

\begin{theorem}\label{the4}
{\sl The set ${\cal T}(X)$ which corresponds to the regular
ultrametric space $X$ with the partial order, defined by inclusion
of balls, is a directed tree where all neighbor vertices are
comparable.

Branching index for vertices of this tree may take finite integer
non--negative values not equal to one, and the maximal vertex (if
exists) has the branching index $\ge 2$. Balls of nonzero diameter
in $X$ correspond to vertices of branching index $\ge 2$ in ${\cal
T}(X)$, and the balls of zero diameter which are maximal subbals
in balls of nonzero diameter correspond to vertices of branching
index 0 in ${\cal T}(X)$. }
\end{theorem}

Remind that a directed set is a partially ordered set, where for
any pair of elements there exists the unique supremum with respect
to the partial order. The ball $I$ the ball in $X$ is a vertex in
${\cal T}(X)$.

Consider the set $X\bigcup {\cal T}(X)$, where we identify the
balls of zero diameter from ${\cal T}(X)$ with the corresponding
points in $X$. We call ${\cal T}(X)$ the tree of balls in $X$, and
$X\bigcup {\cal T}(X)$ the extended tree of balls. One can say
that $X\bigcup {\cal T}(X)$ is the set of all the balls in $X$, of
nonzero and zero diameter.

Introduce the structure of a directed set on $X\bigcup {\cal
T}(X)$. At the tree ${\cal T}(X)$ this structure is the following:
$I<J$ if for the corresponding balls $I\subset J$.

The supremum
$$
{\rm sup}(x,y)=I
$$
of points $x,y\in X$ is the minimal ball $I$, containing the both
points.

Analogously, for $J\in {\cal T}(X)$ and $x\in X$ the supremum
$$
{\rm sup}(x,J)=I
$$
corresponds to the minimal ball $I$, which contains the ball $J$
and the point $x$.

Conversely, starting from a directed tree one can reproduce the
corresponding ultrametric space \cite{Izv}, \cite{nextIzv},
\cite{Lemin}, \cite{Serre}, which will be the absolute of the
directed tree.

Consider a $\sigma$--additive Borel measure $\nu$ with countable
or finite basis on regular ultrametric space $X$, such that for
arbitrary ball $D$ its measure $\nu(D)$ is positive number (i.e.
is not equal to zero).

We study the ultrametric pseudodifferential operator (or the PDO)
of the form considered in \cite{Izv}, \cite{nextIzv}
$$
Tf(x)=\int T{({\rm sup}(x,y))}(f(x)-f(y))d\nu(y)
$$
Here $T{(I)}$ is some nonnegative function on the tree ${\cal T}(
X)$. Thus the structure of this operator is determined by the
direction on $X\bigcup {\cal T}(X)$. This kind of ultrametric PDO
we call {\it the {\rm sup}--operator}.

\subsection{Ultrametric wavelets}

Build a basis in the space $L^2(X,\nu)$ of quadratically
integrable with respect to the measure $\nu$ functions, which we
will call the basis of ultrametric wavelets.

Denote $V_{I}$ the space of functions on the absolute, generated
by characteristic functions of the maximal subballs in the ball
$I$ of nonzero radius. Correspondingly, $V^0_{I}$ is the subspace
of codimension 1 in $V_{I}$ of functions with zero mean with
respect to the measure $\nu$. Spaces $V^0_{I}$ for the different
$I$ are orthogonal. Dimension of the space $V^0_{I}$ is equal to
$p_I-1$.

We introduce in the space $V^0_{I}$ some orthonormal basis
$\{\psi_{Ij}\}$, $j=1,\dots,p-1$. The next theorem shows how to
construct the orthonormal basis in $L^2(X,\nu)$, taking the union
of bases $\{\psi_{Ij}\}$ in spaces $V^0_{I}$ over all non minimal
$I$.

\begin{theorem}\label{basisX}
{\sl 1) Let the ultrametric space $X$ contains an increasing
sequence of embedded balls with infinitely increasing measure.
Then the set of functions $\{\psi_{Ij}\}$, where $I$ runs over all
non minimal vertices of the tree ${\cal T}$, $j=1,\dots,p_I-1$ is
an orthonormal basis in $L^2(X,\nu)$.

2) Let for the ultrametric space $X$ there exists the supremum of
measures of the balls, which is equal to $A$. Then the set of
functions $\{\psi_{Ij}, A^{-{1\over 2}}\}$, where $I$ runs over
all non minimal vertices of the tree ${\cal T}$, $j=1,\dots,p_I-1$
is an orthonormal basis in $L^2(X,\nu)$.
 }
\end{theorem}

The introduced in the present theorem basis we call the basis of
ultrametric wavelets.

\bigskip

The next theorem shows that the basis of ultrametric wavelets is
the basis of eigenvectors for ultrametric PDO.

\begin{theorem}\label{04}{\sl Let the following series converge:
\be\label{seriesconverge} \sum_{J>R} T{(J)} (\nu(J)-\nu({J-1,R}))
<\infty \ee Then the operator
$$
Tf(x)=\int T{({\rm sup}(x,y))}(f(x)-f(y))d\nu(y)
$$
is self--adjoint, has the dense domain in $L^2(X,\nu)$, and is
diagonal in the basis of ultrametric wavelets from the theorem
\ref{basisX}: \be\label{lemma2.1} T\psi_{Ij}(x)=\lambda_I
\psi_{Ij}(x) \ee with the eigenvalues: \be\label{lemma4}
\lambda_{I}=T{(I)} \nu(I)+\sum_{J>I} T{(J)} (\nu(J)-\nu({J-1,I}))
\ee Here $(J-1,I)$ is the maximal vertex which is less than $J$
and larger than $I$.

Also the operator $T$ kills constants. }
\end{theorem}

\subsection{Distributions}

In present Section we introduce the spaces of (complex valued)
test and generalized functions (or distributions) on ultrametric
space $X$. This construction is an analogue of the construction of
the Bruhat--Schwartz space in $p$--adic case, which can be found
in \cite{VVZ}.

\begin{definition}\label{local_constant}{\sl
Function $f$ on (ultrametric) space $X$ is called locally
constant, if for arbitrary point $x\in X$ there exists a positive
number $r$ (depending on $x$), such that the function $f$ is
constant on the ball with center in $x$ and radius $r$:
$$
f(x)=f(y),\qquad \forall y:|xy|\le r
$$
}
\end{definition}

The next definition is an analogue of definition of the
Bruhat--Schwartz space in $p$--adic case.

\begin{definition}\label{test_functions}{\sl
The space of test functions $D(X)$ on ultrametric space $X$ is
defined as the space of locally constant functions with compact
support.}
\end{definition}

Remind that a set (in particular, linear space) is filtrated by a
directed family of subsets (in particular, linear subspaces), if:

1) To any relation of order in this family corresponds inclusion
of subsets. This means that if $A<B$, then exists the embedding of
$A$ into $B$);

2) Any element of the set lies in some set of the family.

In the case of filtrated linear spaces all embeddings are linear.

Consider the tree ${\cal T}={\cal T}(X)$ of balls in ultrametric
space $X$. Introduce the filtration of the space $D(X)$ of test
functions by finite dimensional linear subspaces  $D({\cal S})$,
where ${\cal S}\subset {\cal T}$ is a finite subtree of the
following form.

\begin{definition}\label{wave_type}
{\sl The subset ${\cal S}$ in a directed tree ${\cal T}$ is called
of the regular type, iff:

1) ${\cal S}$ is finite;

2) ${\cal S}$ is a directed subtree in ${\cal T}$ (where the
direction in ${\cal S}$ is the restriction of the direction in
${\cal T}$ onto ${\cal S}$);

3) The directed subtree ${\cal S}$ obey the following property: if
${\cal S}$ contains a vertex $I$ and a vertex $J$: $J<I$,
$|IJ|=1$, then the subtree ${\cal S}$ contains all the vertices
$L$ in ${\cal T}$: $L<I$, $|IL|=1$. }
\end{definition}

Remind that the distance between vertices of a tree is the number
of edges of the path connecting these vertices. The maximal vertex
in ${\cal S}$ we will denote $K$.

We denote by $J$ the ball in $X$, corresponding to vertex $J$ of
${\cal T}(X)$, and by $\chi_J$ we denote the characteristic
function of this ball. In this language we consider a finite set
${\cal S}$ of balls in ultrametric space $X$, characterized by the
properties:

1) If ${\cal S}$ contains balls $I$ and $J$, then it contains
${{\rm sup}(I,J)}$;

2) If ${\cal S}$ contains balls $I$ and $J$: $I\subset J$, then
${\cal S}$ contains all the balls $L$: $I\subset L\subset J$;

3) If ${\cal S}$ contains balls $I$ and $J$, where $J$ is a
maximal subball in $I$, then it contains all the maximal subballs
in $I$.

\begin{definition}\label{D(S)}
{\sl For the finite subtree ${\cal S}\subset {\cal T}$ of the
regular type  consider the space $D({\cal S})$, which is the
linear span of characteristic functions $\chi_J$ with $J\in {\cal
S}$ \footnote{Such subspaces are important in relation to replica
method.}. }
\end{definition}

We consider this space as the subspace in the space $L^2(X,\nu)$
of quadratically integrable with respect to the measure $\nu$
functions on $X$.

Dimension of $D({\cal S})$ is equal to the number of minimal
elements in ${\cal S}$.

The space $D({\cal S})$ has the natural topology (which can be
described as the topology of pointwise convergence), since
$D({\cal S})$ is finite dimensional.  The space $D(X)$ of test
functions on $X$ is the inductive limit of spaces $D({\cal S})$:
\be\label{inductive} D(X)=\lim\,{\rm ind}_{{\cal S}\to {\cal T}}\,
D({\cal S})\ee

Since all spaces $D({\cal S})$ are finite dimensional, the
restriction of the topology of $D({\cal S})$ on any subspace
$D({\cal S}_0)$, ${\cal S}\supset {\cal S}_0$, coincides with the
original topology of $D({\cal S}_0)$. Thus the inductive limit
(\ref{inductive}) is, in fact, the {\it strict inductive limit}
\cite{Schaefer} p.57. By proposition 6.5, \cite{Schaefer}, p.59, a
set $B$ in a strict inductive limit of a countable family of
locally convex spaces $\{E_n\}$ is bounded iff there exists $n$
such that $B\subset E_{n}$ and bounded in it. This implies the
proposition below.

\begin{proposition}\label{topology}{\sl The sequence
$\{f_n\}\in D(X)$ converge, if this sequence lyes in some subspace
$D({\cal S})$ of the filtration and converges (and convergence in
finite dimensional space $D({\cal S})$ is defined uniquely).}
\end{proposition}

We also pay attention that by corollary of theorem 7.4,
\cite{Schaefer}, p.103, we have the following proposition.

\begin{proposition}
{\sl The space $D(X)$ is nuclear.}
\end{proposition}

\begin{definition}\label{generalized_functions}{\sl
Distribution (or generalized function) on $X$ is a linear
functional on the space $D(X)$ of test functions.}
\end{definition}

It is easy to see that this functional automatically will be
continuous (since convergence in the space of test functions is
defined through the convergence in finite dimensional subspaces).
The linear space of generalized functions we denote $D'(X)$. The
convergence in $D'(X)$ is defined as a weak convergence of
functionals. Thus $D'(X)$ is conjugated to $D(X)$ with the weak
topology.

We remind that
$$
D'(X)=\bigcap D'({\cal S})
$$
since $D(X)$ is the strict inductive limit.

\section{Construction of the random field}

Let $X$ be ultrametric space, satisfying the properties 1--3 of
the previous Section (and therefore we have the analysis of
pseudodifferential operators and distribution theory on this
space).

We define the white noise $\phi(x)$ as the generalized Gaussian
real valued random field on $X$ with mean zero and the
$\delta$--like covariation
$$
\langle\phi(x)\phi(x')\rangle=\delta(x-x')
$$
where $\delta$--function is understood in the sense of
distributions in ultrametric space $X$.

Then the Minlos--Sazonov theorem, see the Appendix, implies the
following proposition.

\begin{proposition}
{\sl The Gaussian measure, corresponding to the white noise
$\phi(x)$, is $\sigma$--additive on the $\sigma$--algebra of weak
Borel subsets of the space of distributions $D'(X)$.}
\end{proposition}

Assume we have an ultrametric space $X$, satisfying the necessary
properties, and $T$ is an ultrametric PDO on this space:
$$
Tf(x)=\int T{({\rm sup}(x,y))}(f(x)-f(y))d\nu(y)
$$

Consider the stochastic equation \be\label{SE} T\psi(x)=\phi(x)
\ee where $\phi(x)$ is the white noise.

Remind that the operator $T$ is diagonal in the basis of
ultrametric wavelets $\Psi_{Ij}$:
$$
T\Psi_{Ij}=\lambda_{I}\Psi_{Ij}
$$

Consider the expansion of the white noise over the wavelets
$$
\phi=\sum_{Ij}d_{Ij}\Psi_{Ij}, \qquad d_{Ij}=\int
\phi(x)\Psi_{Ij}(x)d\nu(x)
$$
The coefficients $d_{Ij}$ are mean zero Gaussian random variables,
with the quadratic correlation
$$
\langle d_{Ij}^{*}d_{I'j'}\rangle=\int\int
\langle\phi(x)\phi(x')\rangle
\overline{\Psi_{Ij}(x)}\Psi_{I'j'}(x')
d\nu(x)d\nu(x')=\delta_{II'}\delta_{jj'}
$$
since the wavelets are orthonormal, i.e. $d_{Ij}$ are independent
$\delta$--correlated random variables.

Then by (\ref{SE}) and since the wavelets are basis, we get the
solution \be\label{PW}
\psi(x)=T^{-1}\phi(x)=\sum_{Ij}\lambda^{-1}_{I}d_{Ij}\Psi_{Ij}(x)
\ee

This is a Gaussian mean zero generalized random field (at least
formally).

Compute the quadratic correlation of (\ref{PW}) \be\label{qc}
\langle
\overline{\psi(x)}\psi(x')\rangle=\sum_{Ij}\lambda^{-2}_{I}
\overline{\Psi_{Ij}(x)}\Psi_{Ij}(x') \ee

\begin{lemma}\label{qcl}
{\sl Correlation function (\ref{qc}) has the expression
\be\label{answer} \langle
\overline{\psi(x)}\psi(x')\rangle=-\lambda^{-2}_{{\rm
sup}(x,x')}\nu^{-1}({{\rm sup}(x,x')})+\sum_{I> {\rm
sup}(x,x')}\lambda^{-2}_{I}\left(\nu^{-1}({I-1,{\rm
sup}(x,x')})-\nu^{-1}(I)\right) \ee in the form of the series over
the increasing path in the tree ${\cal T}(X)$, which begins in
${\rm sup}(x,x')$. This expression is correct if and only if both
the series \be\label{conv1} \sum_{J>R} T{(J)}
(\nu(J)-\nu({J-1,R})) \ee and the series \be\label{conv2} \sum_{I>
R}\lambda^{-2}_{I}\left(\nu^{-1}({I-1,R})-\nu^{-1}(I)\right) \ee
where
$$
\lambda_{I}=T{(I)} \nu(I)+\sum_{J>I} T{(J)} (\nu(J)-\nu({J-1,I}))
$$
converge for fixed vertex $R$.

}
\end{lemma}

\noindent{\it Proof}\qquad Taking into account supports of the
wavelets, series (\ref{qc}) over the tree ${\cal T}$ (indexed by
$I$) reduces to the series over the increasing path in the tree,
starting in ${\rm sup}(x,x')$:
$$
\sum_{I,j:I\ge {\rm sup}(x,x')}\lambda^{-2}_{I}
\overline{\Psi_{Ij}(x)}\Psi_{Ij}(x')
$$

Take the sum over $j$. For fixed $I$ let us note that
$$
\sum_{j=0}^{p_I-1}\overline{\Psi_{Ij}}\Psi_{Ij}=\sum_{j=0}^{p_I-1}P_{\chi_{I_j}}
$$
where $P_{\chi_{I}}$ is the orthogonal projector in $L^2(X,\nu)$
onto $\chi_{I}$.

Here $\Psi_{I0}$ is proportional to $\chi_{I}$ and $\Psi_{Ij}$,
$j=0,\dots,p-1$ is the orthonormal basis of wavelets in the space
generated by $\chi_{Ij}$, $j=0,\dots,p-1$ of indicators of the
maximal subbals in $I$.

We have for $I\ge J$
$$
P_{\chi_{I}}\chi_J={\nu(J)\over\nu(I)}\chi_I
$$
or equivalently
$$
P_{\chi_{I}}=\nu^{-1}(I)\overline{\chi_{I}}\chi_{I}
$$

Since
$$
\Psi_{I0}=\nu^{-{1\over2}}(I)\chi_{I}
$$
then
$$
\overline{\Psi_{I0}}\Psi_{I0}=P_{\chi_{I}}
$$

This implies \be\label{sum_in_V_0}
\sum_{j=1}^{p-1}\overline{\Psi_{Ij}}\Psi_{Ij}=\sum_{j=0}^{p_I-1}P_{\chi_{I_j}}-P_{\chi_{I}}
\ee Then, using (\ref{sum_in_V_0}) for sufficiently small balls
$A$ and $B$, containing $x$ and $x'$ correspondingly, we get
$$
\sum_{j=1}^{p-1}\overline{\Psi_{Ij}(x)}\Psi_{Ij}(x')=
\nu^{-1}(A)\nu^{-1}(B)\langle\chi_{A},\left[\sum_{j=0}^{p_I-1}P_{\chi_{I_j}}-P_{\chi_{I}}\right]\chi_{B}\rangle=
$$
$$
=\nu^{-1}({I-1,{\rm sup}(x,x')})\biggr|_{I> {\rm
sup}(x,x')}-\nu^{-1}(I)\biggr|_{I\ge {\rm sup}(x,x')}
$$

This implies the expression (\ref{answer}) for the correlation
function. Investigate, when this expression will be correct (i.e.
when the series will converge).

Eigenvalues of the ultrametric PDO in the basis of wavelets are
$$
\lambda_{I}=T{(I)} \nu(I)+\sum_{J>I} T{(J)} (\nu(J)-\nu({J-1,I}))
$$
Since $\lambda_I$ decreases with the increasing of $I$ (remind
that we consider the case when all $T{(J)}\ge 0$), correlation
function (\ref{answer}) can be divergent (infrared divergence).
The term $\lambda^{-2}_I$ increases with $I$, and $\nu^{-1}(I)$
decreases.

For convergence of the series (\ref{answer}) for all $x$, $x'$ it
is sufficient to have convergence of the series
$$
\sum_{J>R} T{(J)} (\nu(J)-\nu({J-1,R}))
$$
for fixed $R$ (to have existence of all $\lambda_I$), and to have
convergence of the series
$$
\sum_{I>
R}\lambda^{-2}_{I}\left(\nu^{-1}({I-1,R})-\nu^{-1}(I)\right)
$$
for fixed vertex $R$. Moreover, it is easy to see that these
conditions will also be necessary. This finishes the proof of the
lemma.

\begin{lemma}\label{form}{\sl Let conditions (\ref{conv1}) and (\ref{conv2})
be satisfied. Then correlation function (\ref{answer}) of
(\ref{PW}) is a positive continuous bilinear form on $D(X)$. }
\end{lemma}

\noindent{\it Proof}\qquad Since (\ref{qc}) is a sum of orthogonal
projections, multiplied by positive numbers, correlation function
$\langle\overline{\psi(x)}\psi(x')\rangle$ generates positively
definite quadratic form in the space $D(X)$ of test functions in
the case, when the series in (\ref{answer}) converge for all $x$,
$x'$, i.e when conditions (\ref{conv1}) and (\ref{conv2}) are
satisfied. Moreover, this quadratic form will be continuous with
respect to the topology in $D(X)$ (since convergence in $D(X)$
reduces to convergence in finite dimensional subspaces, generated
by finite sets of characteristic functions of balls, and
correlation function (\ref{answer}) is locally constant). This
finishes the proof of the lemma.

\bigskip

Formula (\ref{answer}) gives the analytic expression for the
correlation function of the solution of the stochastic equation
(\ref{SE}). Note that the correlation function (\ref{answer})
depends only on ${\rm sup}(x,x')$, i.e. on the minimal ball in $X$
containing both $x$ and $x'$. This property of independence of
correlations on the details of the process, related to balls which
do not intersect the ball ${\rm sup}(x,x')$, may be discussed as
the ultrametric analogue of the Markov property, see also paper
\cite{BV} for discussion in the $p$--adic case.

The next definition describes how one can introduce Markovianity
for random fields on ultrametric spaces.

\begin{definition}\label{Markov}{\sl
Random field $\psi(x)$ on ultrametric space $X$ with values in
$D'(X)$ is called ultrametric Markovian, if any random variables
$\Psi_I(f)$, $\Psi_J(g)$, where $f,g\in D(X)$ are supported in the
balls $I$ and $J$ correspondingly,
$$
\Psi_I(f)=\int_{I}\psi(x)f(x)d\nu(x), \qquad
\Psi_J(g)=\int_{J}\psi(x)g(x)d\nu(x)
$$
the balls $I$ and $J$ have empty intersection, and at least one of
the functions $f$ and $g$ has zero mean, are independent.

}
\end{definition}

\noindent{\bf Remark}\qquad This is a far generalization of the
Markovianity property for $p$--adic Brownian motion, considered in
\cite{BV}. In \cite{BV} the Markovianity for the random field
$\psi(x)$ of $p$--adic argument was discussed as Markov property
for the sequence $\psi_{\gamma}$, enumerated by integer $\gamma$,
where
$$
\psi_{\gamma}=\int_{|x|_p=p^{\gamma}}\psi(x)dx
$$

\bigskip

Of course, definition \ref{Markov} can also be considered when the
space $X$ is not necessarily ultrametric, but in this general case
we will not have natural examples of such Markovian random fields.

Let us formulate the main result of the present paper.

\begin{theorem}\label{existence}
{\sl Let $X$ be an ultrametric space (satisfying the properties
1--3), $\nu$ a measure on $X$, for which measures of all balls are
positive, $T$ be positive pseudodifferential operator on this
space of the form
$$
Tf(x)=\int T{({\rm sup}(x,y))}(f(x)-f(y))d\nu(y)
$$
and let conditions (\ref{conv1}) and (\ref{conv2}) be satisfied.

Let $\phi(x)$ be white noise on $X$.

Then the solution of the stochastic pseudodifferential equation
$$
T\psi(x)=\phi(x)
$$
is a Gaussian mean zero random field with covariation
$$
\langle \overline{\psi(x)}\psi(x')\rangle=-\lambda^{-2}_{{\rm
sup}(x,x')}\nu^{-1}({{\rm sup}(x,x')})+\sum_{I> {\rm
sup}(x,x')}\lambda^{-2}_{I}\left(\nu^{-1}({I-1,{\rm
sup}(x,x')})-\nu^{-1}(I)\right)
$$

The corresponding Gaussian measure is $\sigma$--additive on the
$\sigma$--algebra of weak Borel subsets of the space of
distributions $D'(X)$. }
\end{theorem}

\noindent{\it Proof}\qquad By lemma \ref{form} we can apply the
Minlos--Sazonov theorem, which provides $\sigma$--additivity of
the Gaussian measure defined by (\ref{qc}) which determines the
random field (\ref{PW}) taking values in $D'(X)$. Lemma \ref{qcl}
gives the expression for the covariation.

This finishes the proof of the theorem.

\begin{theorem}
{\sl Random field $\psi(x)$, defined in theorem \ref{existence} is
ultrametric Markovian in the sense of definition \ref{Markov}. }
\end{theorem}

\noindent{\it Proof}\qquad To prove ultrametric Markovianity we
use Gaussianity of $\psi$ and the formula (\ref{answer}). By
Gaussianity it sufficient to prove Markovianity for the quadratic
correlation function (\ref{answer}).

By (\ref{answer}) the correlation function $\langle
\overline{\psi(x)}\psi(x')\rangle$ is locally constant, and
conditions of definition \ref{Markov} can be proved directly.
Indeed, let us compute
$$
\int_{I}\int_{J}\ov{f(x)}g(x') \langle
\overline{\psi(x)}\psi(x')\rangle d\nu(x)d\nu(x')
$$
where $f$ is supported at the ball $I$, $g$ is supported at the
ball $J$, balls $I$ and $J$ do not intersect and at least one of
the functions $f$ and $g$ has zero mean. Since $\langle
\overline{\psi(x)}\psi(x')\rangle$ is constant for $x\in I$,
$x'\in J$, the above correlation function reduces to
$$
\langle
\overline{\psi(x)}\psi(x')\rangle\int_{I}\ov{f(x)}d\nu(x)\int_{J}g(x')
d\nu(x')=0
$$

This finishes the proof of the theorem.

\bigskip

\noindent{\bf Remark}\qquad Markovianity property for stochastic
processes of real argument, i.e. independence of the future on the
past if the present value of the stochastic process is fixed, is
related to the fact that real numbers are ordered. On $p$--adic
numbers, and, moreover, general ultrametric spaces there is no
natural order. The ultrametric Markovianity in the sense of
definition \ref{Markov} is related to order in the directed tree
of balls ${\cal T}(X)$.

\section{Appendix: the Minlos--Sazonov theorem}

The following material is taken from \cite{DF}.

\begin{definition}
{\sl Let $X$, $Y$ be conjugated linear spaces with the pairing
$\langle\cdot,\cdot\rangle$ and on $X$ there is a quasimeasure
$\mu$. Then the integral
$$
\chi_{\mu}(y)=\int_{X}e^{i\langle x,y\rangle}\mu(dx)
$$
is called the characteristic functional of the quasimeasure $\mu$,
and the map $\mu\mapsto\chi_{\mu}$ is called the Fourier transform
of the quasimeasure.}
\end{definition}

The following is the Minlos--Sazonov theorem.

\begin{theorem}
{\sl Let $X$ be nuclear locally convex space. For the function
$\chi(x)$ to be a characteristic functional of some nonnegative
Radon measure on the space $X'_{\sigma}$ it is sufficient, and if
$X$ is barrelled (tonnele), also is necessary, that the function
$\chi(x)$ is positively defined and continuous in zero in the
topology of the space $X$.}
\end{theorem}

\begin{definition}
{\sl The correlation form of nonnegative cylindric measure $\mu$
with finite second moments
$$
\int_{X'}|\langle y,x\rangle|^2\mu(dy)<\infty, \qquad x\in X
$$
is the bilinear functional on $X$, defined by the expression
$$
B_{\mu}(x_1,x_2)=\int_{X'}\langle x_1,y\rangle\langle
x_2,y\rangle\mu(dy)
$$
}
\end{definition}

\begin{theorem}
{\sl If the correlation form $B_{\mu}(x,y)$ of the cylindric
measure $\mu$ in the space $X'$ is continuously defined in the
nuclear topology $\tau_S(X)$, then $\mu$ is $\sigma$--additive in
$X'$. }
\end{theorem}

\bigskip

\centerline{\bf Acknowledgements}

\medskip

The authors would like to thank I.V.Volovich and A.Kh.Bikulov for
fruitful discussions and valuable comments. This paper has been
partly supported by EU-Network ''Quantum Probability and
Applications''. One of the authors (S.V.Kozyrev) has been partly
supported by The Russian Foundation for Basic Research (project
05-01-00884-a), by the grant of the President of Russian
Federation for the support of scientific schools NSh 1542.2003.1,
by the Program of the Department of Mathematics of Russian Academy
of Science ''Modern problems of theoretical mathematics'', by the
grant of The Swedish Royal Academy of Sciences on collaboration
with scientists of former Soviet Union, and by the grants DFG
Project 436 RUS 113/809/0-1 and RFFI 05-01-04002-NNIO-a.

\end{document}